\documentclass[doublespacing]{elsart}
\usepackage{amsmath}
\usepackage{amsfonts}   % if you want the fonts
\usepackage{amssymb}

\newcommand{\mlabel}[1]{\label{#1}
}

\journal{...}
\newcommand{\seq}{\begin{equation}}                 %start equation
\newcommand{\eeq}[1]{\label{#1}\end{equation}
    }

\newcommand{\epf}{$ \quad \Box$ \par \vspace{1ex}}

\newtheorem{Theorem}{Theorem}[section]
\newcommand{\sthm}{\begin{Theorem}}         %start theorem
\newcommand{\ethm}{\end{Theorem}}           %end theorem

\newtheorem{Corollary}[Theorem]{Corollary}
\newcommand{\scor}{\begin{Corollary}}       %start corollary
\newcommand{\ecor}{\end{Corollary}}         %end corollary

\newtheorem{Lemma}[Theorem]{Lemma}
\newcommand{\slm}{\begin{Lemma}}            %start lemma
\newcommand{\elm}{\end{Lemma}}              %end lemma

\newtheorem{Remark}[Theorem]{Remark}
\newcommand{\srmark}{\begin{Remark}\rm}        %start example
\newcommand{\ermark}{\end{Remark}}             %end example

\newtheorem{Example}[Theorem]{\sc Example}
\newcommand{\sex}{\begin{Example}\rm}        %start example
\newcommand{\eex}{\end{Example}}             %end example

\newcommand{\seql}{\begin{eqnarray*}}       %start a sequence of equations
\newcommand{\eeql}{\end{eqnarray*}}

\newcommand{\smlist}[1]{\begin{list}           %start my list, argum.:\alph,..
                      {(#1{zzcount})}{\usecounter{zzcount}}}
\newcommand{\elist}{\end{list}}

\newcommand{\G}{\Gamma}
\renewcommand{\l}{\lambda}

\newcommand{\s}{\sigma}

\renewcommand{\t}{\tau}
\renewcommand{\O}{\Omega}
\renewcommand{\o}{\omega}

\newcommand{\vect}[1]{\mathbf{#1}}
\providecommand{\norm}[1]{\lVert#1\rVert}
\providecommand{\abs}[1]{\lvert#1\rvert}

\def\e{\varepsilon}

\begin{document}
\begin{frontmatter}

\title{Positive periodic solutions of singular systems of first order ordinary differential equations}
\author{Haiyan Wang}
\ead{wangh@asu.edu}
\address{Division of Mathematical and Natural Sciences\\ Arizona State University \\ Phoenix, AZ 85069-7100, USA}

\pagenumbering{arabic}

\begin{abstract}The existence and multiplicity of positive
periodic solutions for first non-autonomous singular systems are established
with superlinearity or sublinearity assumptions at infinity for an appropriately chosen parameter.
The proof of our results is based on the Krasnoselskii fixed point theorem in a cone.
\end{abstract}


\end{frontmatter}

\section{Introduction}
A variety of population dynamics and physiological processes can be described as the following equation
\begin{equation}\label{eq2}
x'(t)=-a(t)x(t) + \l b(t)f(x(t)).
\end{equation}
Periodic solutions of the type problems have attracted much attention, see e.g. \cite{JIANGWEI2002,OReganWang2005,HWJDE1,WAZLASTA} and
references therein.

On the other hand, recently, there are a considerable interest in the existence of positive periodic solutions of singular systems
of the second order differential
equations, see Chu, Torres and Zhang \cite{chujde2007}, Franco and Webb \cite{Franco2006}, Jiang, Chu, O'Regan, and Agarwal \cite{jcoa}, and the author \cite{Wang2010} and references therein.  It has been shown that many results for nonsingular systems still valid for
singular cases. In particular, the author \cite{Wang2010} demonstrates that the Krasnoselskii fixed point theorem on compression and expansion of cones
can be effectively used to deal with  singular problems. In fact, by choosing appropriate cones, the singularity of the systems is essentially removed and the associated operator becomes well-defined for
certain ranges of functions even there are negative terms.

Agarwal and O'Regan  \cite{Agarwal2003} provided some results on solutions of singular first order differential equations.
Chu and  Nieto \cite{chu2008} showed
the existence of periodic solutions for singular first order differential equations with impulses based on a nonlinear alternative of Leray-Schauder. 
The results in \cite{Agarwal2003,chu2008} for first order differential equations deal with a single equation. 
In this paper, by employing the Krasnoselskii fixed point theorem on compression and expansion of a cone,
we shall establish the existence and multiplicity of positive $\o$-periodic solutions for the following singular non-autonomous $n$-dimensional system
\begin{equation}\label{eq1}
x_i(t)=-a_i(t)x_i(t) + \l b_i(t)f_i(x_1(t),..., x_n(t)), i=1,...,n.
\end{equation}
where $\l>0$ is a positive parameter. Our results give an almost complete structure of the existence
of positive periodic  solutions of (\ref{eq1}) with an appropriately chosen parameter. Our results further show that
there are analogous results between the first order and second ordinary differential equations.

First we make assumptions for (\ref{eq1}). Let $\mathbb{R}=(-\infty, \infty)$, $\mathbb{R}_+=[0, \infty)$,
$\mathbb{R}_+^n=\Pi_{i=1}^n \mathbb{R}_+$ and for any
$\vect{u}=(u_1,...,u_n) \in \mathbb{R}^n_+$,
$\norm{\vect{u}}=\sum_{i=1}^n \abs{u_i}$.

{(H1)~~$a_i, b_i$ $\in C(\mathbb{R}, [0, \infty))$ are
$\o$-periodic functions such that $\int_0^{\o}a_i(t)dt >0$,
$\int_0^{\o}b_i(t)dt >0$, $i=1,\dots,n$.}

{(H2)~~$f_i: \mathbb{R}_+^n \setminus \{0\} \to (0,\infty)$ is continuous, $i=1,\dots,n$}\\

Our main results are:

\sthm\mlabel{th1} Let (H1),(H2) hold. Assume that $\lim_{ \norm{\vect{u}} \to 0} f_i(x)=\infty$ for some $ i=1,...,n.$ \\
(a). If $\lim_{ \norm{\vect{u}} \to \infty} \frac{f_i(x)}{\norm{\vect{u}}}=0$,  $i=1,\dots,n$ , then, for all $\l>0$, (\ref{eq1}) has
a positive periodic solution. \\
(b). If $\lim_{ \norm{\vect{u}} \to \infty} \frac{f_i(x)}{\norm{\vect{u}}}=\infty$  for $i=1,\dots,n$, then, for all sufficiently small $\l>0$, (\ref{eq1}) has
two positive periodic solutions. \\
(c). There exists a $\l_0>0$ such that (\ref{eq1}) has a positive periodic solution for $0 < \l < \l_0$.
\ethm

\srmark\label{rem1}
As discussed in \cite{Wang2010}, we can extend Theorem \ref{th1} to the following singular systems with possible negative $e_i$,
\begin{equation}\label{eq11}
x_i'(t)=-a_i(t)x_i(t) + \l b_i(t)f_i(x_1(t),..., x_n(t)) + \l e_i(t), i=1,...,n.
\end{equation}
where $e_i(t), i=1,...,n,$ are continuous $\o$-periodic functions. When $e_i(t)$ takes negative values, we shall need
a stronger condition on $b_i (b_i>0).$
%\sthm\mlabel{th5} Let \rm{(H1),(H2) and $b_i(t)>0$ for $t \in [0, T]$, $i=1,\dots,n.$} Assume that $\lim_{ \norm{\vect{u}} \to 0} f_i(x)=\infty$ for
%$i=1,...,n.$ \\
%(a). If $\lim_{ \norm{\vect{u}} \to \infty} f_i(x)=\infty$ and  $\lim_{\norm{\vect{u}}\to \infty} \frac{f_i(x)}{\norm{\vect{u}}}=0$,  $i=1,\dots,n$ ,
% then there exists  $\l_0>0$ such that
%(\ref{eq11}) has a positive periodic solution for $\l > \l_0$. \\
%(b). If $\lim_{\norm{\vect{u}} \to \infty} \frac{f_i(x)}{\norm{\vect{u}}}=\infty$  for  $i=1,\dots,n$ , then, for all sufficiently small $\l>0$, (\ref{eq11})  has
%two positive periodic solutions.\\
%(c). There exists a $\l_1>0$ such that (\ref{eq11})  has a positive periodic solution for $0 < \l < \l_1$. \\
%\ethm

Such a result  can be proved in the same way as in \cite{Wang2010}. We will not give a detailed proof here.  The idea to deal with negative $e_i$ is to split $b_i(s) f_i(x(s))+e_i(t)$ into the two terms $\frac{1}{2}b_i(s) f_i(x(s))$
and $\frac{1}{2}b_i(s) f_i(x(s))+e_i(t)$.
The first term is always nonnegative and used to carry out the estimates of the operator. We will make the second term
$\frac{1}{2}b_i(s) f_i(x(s))+e_i(t)$ nonnegative by choosing appropriate domains of $f_i$. This is possible because
$\lim_{ \norm{\vect{u}} \to 0} f_i(x)=\infty$ or $\lim_{ \norm{\vect{u}} \to \infty} f_i(x)=\infty$.
The choice of the even split
of $b_i(s) f_i(x(s))$ here is not necessarily optimal in terms of obtaining maximal $\l$-intervals for the existence of periodic solutions of the systems.
\ermark

\srmark\label{rem2}
O'Regan and the author \cite{OReganWang2005}, and the author \cite{HWJDE1} established the existence, multiplicity and nonexistence of positive
periodic solution of the first order ODE
\begin{equation}\label{eq90}
x_i'(t)=a_i(t)g_i(x(t))x_i(t) - \l b_i(t)f_i(x(t-\t(t))), i=1,...,n.
\end{equation}
where $g_i$ are positive bounded functions and $\t \in C(\mathbb{R}, [0, \infty))$ is a $\o$-periodic function. These results can also be extended to (\ref{eq90})
if $f_i$ has a singularity at zero.
\ermark

%%%%%%%%%%%%%%%%%%%%%%%%%%%%%%%%%%%%%%%%%%%%%%%%%%%%%%%%%%%

\setcounter{equation}{0}
\section{Preliminaries}

We recall some concepts and conclusions of an operator in a cone. Let $E$
be a Banach space and $K$ be a closed, nonempty subset of $E$. $K$
is said to be a cone if $(i)$~$\alpha u+\beta v\in K $ for all
$u,v\in K$ and all $\alpha,\beta \geq 0$ and $(ii)$~$u,-u\in K$ imply
$u=0$. The following well-known result of the fixed
point theorem is crucial in our arguments.

\slm\mlabel{lm1} {\rm (\cite{KRAS})} Let $X$ be a Banach
space and $K\ (\subset X)$ be a cone. Assume that $\Omega_1,\
\Omega_2$ are bounded open subsets of $X$ with $0 \in \Omega_1,\bar\Omega_1 \subset \Omega_2$, and let
$$
\mathcal{T}: K \cap (\bar{\Omega}_2\setminus \Omega_1 ) \rightarrow K
$$
be completely continuous operator such that either
\begin{itemize}
\item[{\rm (i)}] $\| \mathcal{T}u \| \geq \| u \|,\ u\in K\cap \partial
     \Omega_1$ and $ \| \mathcal{T}u \| \leq \| u \|,\ u\in K\cap \partial
     \Omega_2$; or

\item[{\rm (ii)}] $\| \mathcal{T}u \| \leq \| u \|,\ u\in K\cap \partial
     \Omega_1$ and $\| \mathcal{T}u \| \geq \| u \|,\ u\in K\cap \partial
     \Omega_2$.
\end{itemize}
Then $\mathcal{T}$ has a fixed point in $K \cap ( \bar \Omega_2 \backslash
     \Omega_1)$.

\elm

We now introduce some notation. For $r>0$, let
$$\begin{array}{l}
\s=\min\limits_{i=1,\dots,n}\{\s_i \}>0,~~{\rm where}~~\s_i =
e^{-\int_0^{\o}a_i(t)dt}, i=1,\dots,n,\\
M(r)=\max\{f_i(\vect{u}): \vect{u} \in
\mathbb{R}_+^n, \s r\leq \norm{\vect{u}} \leq r,i=1,\dots,n\}>0,\\
m(r)=\min\{f_i(\vect{u}): \vect{u} \in
\mathbb{R}_+^n, \s r\leq  \norm{\vect{u}} \leq r,i=1,\dots,n\}>0,\\
\G=\s \min\limits_{n=1,\dots,n}\{\frac{\int^{\o}_{0} b_i(s)
ds}{\s^{-1}_i-1}\}>0,~~\chi = \sum_{i=1}^n
\frac{\s^{-1}_i}{\s^{-1}_i-1} \int^{\o}_{0} b_i(s)ds>0.
\end{array}
$$

In order to apply Lemma \ref{lm1} to (\ref{eq1}), let $X$ be the
Banach space defined by
$$X=\{\vect{u}(t)\in C(\mathbb{R},\mathbb{R}^n):
\vect{u}(t+\o)=\vect{u}(t), t \in \mathbb{R}, i=1,\dots,n\}$$ with
a norm $\displaystyle{\norm{\vect{u}}= \sum_{i=1}^n
\sup_{t\in[0,\o]} \abs{u_i(t)}},$ for $\vect{u}=(u_1,...,u_n) \in
X.$  For $\vect{u} \in X$ or $\mathbb{R}^n_+$, $\norm{\vect{u}}$
denotes the norm of $\vect{u}$ in $X$ or $\mathbb{R}^n_+$,
respectively.

Define
$$
K = \{\vect{u}=(u_1,...,u_n) \in X: u_i(t) \geq \s_i
\sup_{t\in[0,\o]} \abs{u_i(t)}, i=1,\dots,n, t \in [0, \o] \}.
$$
It is clear $K$ is cone in $X$ and $\min_{ t \in [0,\o]} \sum_{i=1}^n |u_i(t)| \geq \s \norm{u}$ for $\vect{u}=(u_1,...,u_n) \in K$.
For $r>0$, define $\O_r =
\{\vect{u} \in K: \norm{\vect{u}} < r \}. $ It is clear that
$\partial \O_r = \{\vect{u} \in K: \norm{\vect{u}}=r\}$.  Let
$\vect{T}_{\l}: K \setminus \{0\}\to X$ be a map with components
$(T_{\l}^1,...,T_{\l}^n)$:
\begin{equation}\label{T_def}
T_{\l}^i\vect{u}(t) = \l \int^{t+\o}_t G_i(t, s) b_i(s)
f_i(\vect{u}(s))ds,~~i=1,\dots,n,
\end{equation}
where
$$
G_i(t,s)=\frac{e^{\int_t^s a_i(\theta)d\theta}}{\s^{-1}_i-1}
$$
satisfying
$$
\frac{1}{\s^{-1}_i-1} \leq G_i(t,s) \leq
\frac{\s^{-1}_i}{\s^{-1}_i-1}\; \text{for} \; t \leq s \leq t+\o.
$$
\slm\mlabel{lm-compact} Assume \rm{(H1)-(H2)} hold. Then $\vect{T}
_{\l}(K \setminus \{0\}) \subset K$ and $\vect{T}_{\l}: K\setminus \{0\} \to K$ is compact and
continuous. \elm

\pf If $ \vect{u}=(u_1,...,u_n) \in K \setminus \{0\}$, then $\min_{ t \in [0,\o]} \sum_{i=1}^n |u_i(t)| \geq \s \norm{\vect{u}}>0$,
and then $T _{\l}^i$ is defined. In view of the definition of $K$, for
$\vect{u} \in K \setminus \{0\}$, we have, $i=1,\dots,n,$
\begin{equation*}
\begin{split}
(T_{\l}^i\vect{u})(t+\o)
& =\l \int^{t+2\o}_{t+\o} G_i(t+\o, s) b_i(s) f_i(\vect{u}(s))ds \\
& =\l \int^{t+\o}_{t} G_i(t, s) b_i(s) f_i(\vect{u}(s))ds
=(T_{\l}^i\vect{u})(t).
\end{split}
\end{equation*}
It is easy to see that $\int^{t+\o}_{t} b_i(s)
f_i(\vect{u}(s))ds$ is a constant because of the periodicity
of $b_i(t) f_i(\vect{u}(t))$. One can show that, for
$\vect{u} \in K \setminus \{0\} $ and $t \in [0,\o]$, $i=1,\dots,n,$
\begin{equation*}
\begin{split}
T_{\l}^i\vect{u}(t)
& \geq \frac{1}{\s^{-1}_i-1}\l \int^{t+\o}_{t} b_i(s) f_i(\vect{u}(s))ds \\
& = \s_i \frac{\s^{-1}_i}{\s^{-1}_i-1} \l \int^{\o}_{0} b_i(s)
f_i(\vect{u}(s))ds\geq \s_i  \sup_{t\in[0,\o]}
\abs{T_{\l}^i\vect{u}(t)}.
\end{split}
\end{equation*}
Thus $\vect{T} _{\l}(K \setminus \{0\}) \subset K$ and it is easy to show that
$\vect{T}_{\l}:  K \setminus \{0\}\to K$ is compact and continuous. \epf

\slm\mlabel{lm-fixed-equation-equal} Assume that \rm{(H1)-(H2)}
hold. Then $\vect{u}\in K \setminus \{0\} $ is a positive periodic solution of
(\ref{eq1}) if and only if it is a fixed point of $\vect{T} _{\l}$
in $K \setminus \{0\}$

\elm \pf If $\vect{u}=(u_1,\dots,u_n) \in K \setminus \{0\}$ and
$\vect{T}_{\l}\vect{u}=\vect{u}$, then, for $i=1,\dots,n,$
\begin{equation*}
\begin{split}
u_i'(t)
& = \frac{d}{dt} (\l \int^{t+\o}_t G_i(t, s) b_i(s) f_i(\vect{u}(s))ds)\\
& =  \l G_i(t, t+\o) b_i(t+\o)f_i(\vect{u}(t+\o)- \l G_i(t,t)b_i(t)f_i(\vect{u}(t)) \\
& \quad - a_i(t) T_{\l}^iu(t) \\
& =  \l [G_i(t, t+\o)-G_i(t,t)]b_i(t)f_i(\vect{u}(t))-a_i(t) T_{\l}^iu(t)\\
& = - a_i(t) u_i(t)+ \l b_i(t)f_i(\vect{u}(t)).
\end{split}
\end{equation*}
Thus $\vect{u}$ is a positive $\o$-periodic solution of
(\ref{eq1}). On the other hand, if $\vect{u}=(u_1,\dots,u_n)$ is a
positive $\o$-periodic function, then $ \l
b_i(t)f_i(\vect{u}(t))= a_i(t) u_i(t)+u_i'(t)$ and
\begin{equation*}
\begin{split}
T_{\l}^i\vect{u}(t)
& = \l \int^{t+\o}_t G_i(t, s) b_i(s) f_i(\vect{u}(s))ds\\
& =  \int^{t+\o}_t G_i(t, s) (a_i(s) u_i(s)+u_i'(s))ds \\
& =  \int^{t+\o}_t G_i(t, s) a_i(s) u(s)ds +\int^{t+\o}_t G_i(t, s) u'_i(s)ds \\
& =  \int^{t+\o}_t G_i(t, s) a_i(s) u(s)ds + G_i(t, s) u(s)|^{t+\o}_t - \int^{t+\o}_t G_i(t, s) a_i(s) u_i(s)ds \\
& = u_i(t).
\end{split}
\end{equation*}
Thus, $\vect{T}_{\l}\vect{u}=\vect{u}$, Furthermore, in view of
the proof of Lemma \ref{lm-compact}, we also have $u_i(t) \geq
\s_i \sup_{t\in[0,\o]} u_i(t)$ for $t \in [0,\o].$ That is,
$\vect{u}$ is a fixed point of $\vect{T}_{\l}$ in $K \setminus \{0\}$. \epf

\slm\mlabel{f_estimate_>*} Assume that \rm{(H1)-(H2)} hold. For
any $\eta > 0$ and $ \vect{u}=(u_1, \dots, u_n) \in K \setminus \{0\} $, if there
exists a $f_i$ such that
\mbox{$f_i(\vect{u}(t)) \geq \sum_{j=1}^n u_j(t)\eta$ } for $ t
\in [0, \o]$, then $ \norm{\vect{T}_{\l}\vect{u}} \geq \l \G
\eta \norm{\vect{u}}. $ \elm

 \pf Since $\vect{u} \in K\setminus \{0\}$ and
\mbox{$f_i(\vect{u}(t)) \geq \sum_{j=1}^n u_j(t)\eta$ } for $ t
\in [0, \o]$, we have
\begin{equation*}
\begin{split}
(T_{\l}^i\vect{u})(t)
& \geq \frac{1}{\s^{-1}_i-1}\l \int^{\o}_{0} b_i(s) f_i(\vect{u}(s))ds \\
& \geq \frac{1}{\s^{-1}_i-1}\l \int^{\o}_{0} b_i(s) \sum_{j=1}^n u_j(s) \eta ds \\
& \geq \frac{1}{\s^{-1}_i-1}\l \int^{\o}_{0} b_i(s) ds \sum_{j=1}^n \s_j\sup_{t\in[0,\o]} u_j(t) \eta  \\
& \geq \l \min_{i=1,\dots,n}\{\s_i\} \frac{\int^{\o}_{0} b_i(s)
ds}{\s^{-1}_i-1} \eta  \norm{\vect{u}}.
\end{split}
\end{equation*}
Thus $\norm{\vect{T}_{\l}\vect{u}} \geq \l \G \eta
\norm{\vect{u}}$. \epf

Let $\hat{f}_i: [1, \infty) \to \mathbb{R}_+$ be the function given by
$$\hat{f}_i(\theta) =\max \{f_i(u):u \in \mathbb{R}_+^n \; \rm{and}\; 1 \leq \abs{u} \leq \theta \}, i=1,...,n.$$
It is easy to see that $\hat{f}_i(\theta)$ is a nondecreasing function on $[1,\infty)$. The following lemma is essentially the same  as \cite[Lemma 3.6]{Wang2010} and \cite[Lemma 2.8]{HWJMAA1}.
\slm (\cite{Wang2010,HWJMAA1})\mlabel{lm6}
Assume \rm{(H1)} holds. If $\lim_{ \abs{x} \to \infty} \frac{f_i(x)}{\abs{x}}$ exists (which can be infinity), then $\lim_{\theta \to \infty} \frac{\hat{f}_i(\theta)}{\theta}$ exists and
$\lim_{\theta \to \infty} \frac{\hat{f}_i(\theta)}{\theta}=\lim_{ \abs{x} \to \infty} \frac{f_i(x)}{\abs{x}}$.
\elm

\slm\mlabel{f_estimate_<*} Assume that
\rm{(H1)-(H2)} hold.  Let $r > \frac{1}{\s}$ and if there exists an $\e > 0$ such that
$$
\hat{f}_i(r) \leq  \e r,\;\; i=1,...,n,
$$
then
$$
\norm{\vect{T}_{\l}\vect{u}} \leq \l \chi \varepsilon
\norm{\vect{u}} $$
for $\vect{u}=(u_1, \dots, u_n) \in \partial\O_{r}$.
\elm

\pf From the definition of $\vect{T}$, for
$\vect{u} \in \partial\O_{r}$, we have
\begin{equation*}
\begin{split}
\norm{\vect{T}_{\l}\vect{u}}
  & \leq  \sum_{i=1}^n \frac{\s^{-1}_i}{\s^{-1}_i-1} \l \int^{\o}_{0} b_i(s) f_i(\vect{u}(s))ds  \\
  & \leq  \sum_{i=1}^n \frac{\s^{-1}_i}{\s^{-1}_i-1} \l \int^{\o}_{0} b_i(s) \hat{f}_i(r)ds  \\
  & \leq   \sum_{i=1}^n \frac{\s^{-1}_i}{\s^{-1}_i-1} \l \int^{\o}_{0} b_i(s)ds  \varepsilon  \norm{\vect{u}}=  \l \chi \varepsilon \norm{\vect{u}}.
\end{split}
\end{equation*}
\epf

In view of the definitions of $m(r)$ and $M(r)$,  it follows that
$ M(r) \geq f_i(\vect{u}(t)) \geq  m(r)$ \; $\rm{for}\; t \in [0,
\o]$, $i=1,\dots,n$ if $\vect{u} \in \partial \O_{r}$, $r>0$ .
Thus it is easy to see that the following two lemmas can be shown
in similar manners as in Lemmas \ref{f_estimate_>*} and
\ref{f_estimate_<*}.

\slm\mlabel{lm8} Assume \rm{(H1)-(H2)} hold. If $ \vect{u} \in
\partial \O_{r}$, $r
>0$, then $ \norm{\vect{T}_{\l}\vect{u}}   \geq \l \frac{\G}{\s} m(r). $\elm

%\pf Since $ f_i(\vect{u}(t)) \geq  m(r) \; \rm{for}\; t \in
%[0, \o],$ $i=1,\dots,n,$ it is easy to see that this lemma can be
%shown in a similar manner as in Lemma \ref{f_estimate_>*}. \epf

\slm\mlabel{lm9} Assume \rm{(H1)-(H2)} hold. If $ \vect{u}\in
\partial \O_{r}$, $r >0$, then $ \norm{\vect{T}_{\l}\vect{u}}
\leq \l \chi M(r). $ \elm

%\pf Since $ f_i(\vect{u}(t)) \leq M(r)
%\; \rm{for}\; t \in [0,\o]$, $i=1,\dots,n,$ it is easy to see that
%this lemma can be shown in a similar manner as in Lemma
%\ref{f_estimate_<*}. \epf

\setcounter{equation}{0}
\section{Proof of Theorem \ref{th1}}

Part (a). From the assumptions,
there is an $r_1 > 0$ such that
$$
f_i(\vect{u}) \geq  \eta \norm{\vect{u}}
$$
for $ \vect{u}=(u_1,...,u_n) \in \mathbb{R}_+^n$ and $ 0<\norm{\vect{u}} \leq r_1,$
where $\eta > 0$ is chosen so that
$$
\l \G \eta > 1.
$$
If $ \vect{u}=(u_1,...,u_n) \in  \partial \O_{r_1}$, then
$$f_i(\vect{u}(t)) \geq \eta\sum_{i=1}^n u_i(t), \;\; {\rm for } \;\; t \in [0,1].$$
Lemma ~\ref{f_estimate_>*} implies that
$$
\norm{\vect{T}_{\l}\vect{u}} \geq \l \G \eta \norm{\vect{u}} > \norm{\vect{u}} \quad \textrm{for} \quad  \vect{u} \in \partial\O_{r_1}.
$$
We now determine $\O_{r_2}$.  Since $\lim_{ \norm{\vect{u}} \to \infty} \frac{f_i(x)}{\norm{\vect{u}}}=0$,  $i=1,\dots,n$
it follows from Lemma ~\ref{lm6} that $\lim_{\theta \to \infty }\frac{\hat{f}_i(\theta)}{\theta}=0$, $i=1,...,n.$
Therefore there is an $r_2>\max\{2r_1,\frac{1}{\s}\}$ such that
$$
\hat{f}_i(r_2) \le \e r_2,\;i=1,...,n,
$$
where the constant $\e > 0$ satisfies
$$
\l \e \chi < 1.
$$
Thus, we have by Lemma ~\ref{f_estimate_<*} that
$$
\norm{\vect{T}_{\l}\vect{u}} \leq \l \e \chi \norm{\vect{u}} < \norm{\vect{u}} \quad \textrm{for} \quad  \vect{u} \in \partial\O_{r_2}.
$$
By Lemma ~\ref{lm1},
It follows that
$\vect{T}_{\l}$ has a fixed point in  $\O_{r_2} \setminus \bar{\O}_{r_1}$, which is the desired positive solution of (\ref{eq1}).
\epf

Part (b).
Fix a number $r_1 > 0$. Lemma \ref{lm9} implies that there exists a $\l_0 >0$ such that
$$
\norm{\vect{T}_{\l}\vect{u}}  <  \norm{\vect{u}}, \; {\rm for} \; \vect{u} \in  \partial \O_{r_1},\; 0< \l < \l_0.
$$
In view of  $\lim_{ \norm{\vect{u}} \to 0}f_i(x) = \infty$, there is a positive number $r_2 <  r_1$ such that
$$
f_i(\vect{u}) \geq \eta \norm{\vect{u}}
$$
for $ \vect{u}=(u_1,...,u_n) \in \mathbb{R}_+^n$ and $ 0<\norm{\vect{u}} \leq r_2,$
where $\eta > 0$ is chosen so that
$$
\l \G \eta > 1.
$$
Then
$$f_i(\vect{u}(t)) \geq \eta\sum_{i=1}^n u_i(t),$$
for $\vect{u}=(u_1,...,u_n) \in  \partial \O_{r_2}, \;\;t \in [0,1].$
Lemma ~\ref{f_estimate_>*} implies that
$$
\norm{\vect{T}_{\l}\vect{u}} \geq \l \G \eta \norm{\vect{u}} > \norm{\vect{u}} \quad \textrm{for} \quad  \vect{u} \in \partial\O_{r_2}.
$$
On the other hand, since $\lim_{\norm{\vect{u}} \to \infty} \frac{f_i}{\norm{\vect{u}}}=\infty$, there is
an $\hat{H} > 0$ such that
$$
f_i(\vect{u}) \geq \eta \norm{\vect{u}}
$$
for $ \vect{u}=(u_1,...,u_n) \in \mathbb{R}_+^n$ and $\norm{\vect{u}} \geq \hat{H}$ ,
where $\eta > 0$ is chosen so that
$$
\l \G \eta > 1.
$$
Let $r_3 = \max\{2r_1, \frac{\hat{H}}{\s} \}$. If $ \vect{u}=(u_1 ,...,u_n) \in \partial \O_{r_3}$, then
$$ \min_{0\leq t \leq\o} \sum_{i=1}^n u_i(t) \geq \s
\norm{\vect{u}}= \s r_3 \geq \hat{H},$$
which implies that
$$
f_i(\vect{u}(t)) \geq  \eta\sum_{i=1}^n u_i(t)\; \rm{for} \; t \in [0,\o].
$$
It follows from Lemma ~\ref{f_estimate_>*} that
$$
\norm{\vect{T}_{\l}\vect{u}} \geq \l\G \eta \norm{\vect{u}} > \norm{\vect{u}} \quad \rm{for} \quad  \vect{u} \in \partial\O_{r_3}.
$$
It follows from Lemma ~\ref{lm1} that
$\vect{T}_{\l}$ has two fixed points $\vect{u}_1 $ in  $\O_{r_1} \setminus \bar{\O}_{r_2}$ and $\vect{u}_2 $ $ \in \O_{r_3} \setminus \bar{\O}_{r_1}$ such that
$$
r_2 < \norm{\vect{u}_1} < r_1 < \norm{\vect{u}_2} < r_3.
$$
Consequently, (\ref{eq1})
has two positive solutions for $ 0< \l < \l_0$.

Part (c).
Fix a number $r_1 > 0$. Lemma \ref{lm8} implies that there exists a $\l_0 >0$ such that
$$
\norm{\vect{T}_{\l}\vect{u}}  <  \norm{\vect{u}}, \; {\rm for} \; \vect{u} \in  \partial \O_{r_1},\; 0< \l < \l_0.
$$
In view of  $\lim_{x\to 0}f_i(x) = \infty$, there is a positive number $r_2 <  r_1$ such that
$$
f_i(\vect{u}) \geq \eta \norm{\vect{u}}
$$
for $ \vect{u}=(u_1,...,u_n) \in \mathbb{R}_+^n$ and $ 0<\norm{\vect{u}} \leq r_2,$
where $\eta > 0$ is chosen so that
$$
\l \G \eta > 1.
$$
Then
$$f_i(\vect{u}(t)) \geq \eta\sum_{i=1}^n u_i(t),$$
for $\vect{u}=(u_1,...,u_n) \in  \partial \O_{r_2}, \;\;t \in [0,1].$
Lemma ~\ref{f_estimate_>*} implies that
$$
\norm{\vect{T}_{\l}\vect{u}} \geq \l \G \eta \norm{\vect{u}} > \norm{\vect{u}} \quad \textrm{for} \quad  \vect{u} \in \partial\O_{r_2}.
$$
It follows from Lemma ~\ref{lm1} that
$\vect{T}_{\l}$ has a fixed point in  $\O_{r_1} \setminus \bar{\O}_{r_2}$. Consequently, (\ref{eq1})
has a positive solution for $ 0< \l < \l_0$.

\end{document}